\documentclass[a4paper]{article}
\usepackage{a4wide}
\usepackage{amsmath}	
\usepackage{amsfonts}
\usepackage{amssymb}
\usepackage{amsthm}

\newtheorem{thm}{Theorem}
\newtheorem{lem}[thm]{Lemma}
\newenvironment{pf}{\noindent{\bf Proof: }}{\qed\medskip\par}

\def\eps{\varepsilon}

\def\R{{\mathbb R}}
\def\Z{{\mathbb Z}}
\def\grad{\bigtriangledown\!}

\begin{document}
\title{Families of surfaces lying in a null set}
\author{Laura Wisewell}
\maketitle
\section{Introduction}
In this note we generalise the following result of Sawyer
\cite{sawyer:planecurves}:
\begin{thm}
\label{saw}
There is a function $\psi$ on $\R$ such that whenever $g$ is a
real-valued Borel measurable function on (a subset of)
$\R\times\R^{n-1}$ with the property that $y\mapsto g(y,t)$ is $C^1$
for a.e.~$t$, the set
 $$E_f:=\bigcup_y\left\{(x,t)\in\R\times\R^{n-1}:x=g(y,t)-\psi(y)\right\}$$ has measure zero.
\end{thm}
That is, a smooth one-parameter family of measurable hypersurfaces may
be translated to lie in a null set.  Moreover, the translations may be
taken parallel to $\R$ and need not depend on $g$.

Our motivation came from studying curved analogues of the Kakeya and Nikodym problems: we wished to show that a null set could contain a translate of every member of a specified family of curves, or in the Nikodym case, a curve from a specified family through every point.
So we generalise Sawyer's result in two ways.
First, we will allow any codimension
since curves of course have codimension $n-1$.  Second, we do not want
to be restricted to using translations, since in the Nikodym case it
is the ``directions'' or ``shapes'' we are allowed to vary while the positions are kept
fixed.  So we remove all distinction between ``shape parameters'' and
``position parameters'', simply denoting those that are ``given'' by
$y$ and those we are free to choose by $\omega$.

Consider  objects of the following form
$$\Gamma(y,\omega):=\left\{\begin{pmatrix}f(y,\omega,t)\\t\end{pmatrix}:t\in\R^d\right\}$$
where $f:\R^p\times\R^q\times\R^d\to\R^{n-d}$.  So $\Gamma (y,\omega)$
can be thought of as a $d$-dimensional surface in $\R^n$, and the family of them has
$p+q$ parameters in total.  In the case above,
$f(y,\omega,t)=g(y,t)-\omega$.

Our aim is to show that under certain hypotheses, a null set may
include a representative of every combination of the first $p$
parameters provided that the remaining $q$ parameters can be chosen to
depend on them.  That is, there exists a set of measure zero that
includes a $\Gamma (y,\omega(y))$ for every $y$.  In fact, this
function of $y$ will be the obvious generalisation of Sawyer's
universal translation function $\psi$, and will not depend on $f$.

More precisely, our theorem is the following
\begin{thm}
\label{thm:null} There is a function $\psi:\R^p\to\R^q$ with the following property: 
 Let $f:\R^p\times\R^q\times\R^d\to\R^{n-d}$ where $p\leq n-d\leq q$
 and $d<n$.  Suppose that $f$ is measurable, and for almost every fixed $t$ that the
 map $(y,\omega)\mapsto f(y,\omega,t)$ is $C^1$, that the Jacobian
 $\frac{\partial f}{\partial\omega}$ has full rank (namely
 $n-d$) and that this Jacobian is Lipschitz.  Then the set
 $$E_f:=\bigcup_{y\in\R^p}\Gamma (y,\psi(y))$$
 has measure zero.
\end{thm}
Thus Theorem~\ref{saw} is the special case obtained by setting $d=n-1$, $p=q=1$ and $f(y,\omega,t)=g(y,t)-\omega$.

The proof will have three parts.  First, we define the universal
 transformation function $\psi$, which will be the obvious higher
 dimensional analogue of that used by Sawyer.  Then, we show that all
 of the slices through the set at fixed $t$ have zero measure, which
 is where the conditions on the Jacobian and the $C^1$ assumption are used.  Finally we
 show that the whole set is measurable, using the $C^1$ condition again.  This allows us to apply Fubini's theorem to
 obtain the result.
\section[Definition of $\psi$]{Definition of $\boldsymbol{\psi}$}
 We begin with a few easily verified facts needed for the proof.
\begin{enumerate}
\item Factorial Expansion: Every $a\in(0,1]$ has a unique expansion of the form
 $$a=\sum_{n=2}^\infty\frac{a_n}{n!}$$ where the $a_n$ are integers,
$0\leq a_n\leq n-1$ and infinitely many of the coefficients are
non-zero.
\item There are countably many numbers in $(0,1]$ that also have a finite factorial expansion. 
\item ${\displaystyle \sum_N^\infty\frac{n-1}{n!} = \frac{1}{(N-1)!}}$
\end{enumerate}
 All norms, whether of matrices or vectors, will denote the largest
 absolute value of the entries---this is merely to avoid keeping
 track of constants, since of course all norms on a finite dimensional
 space are equivalent.

 We shall use subscripts to denote the coefficients of the factorial
 expansions of the vectors $y$ rather than their components.  Thus
 for $y\in(0,1]^p$ we can write
 $y=\sum_{n=2}^\infty\frac{y_n}{n!}$ in the natural way.

 Our aim is to construct a kind of ``universal transformation
 function'' $\psi:\R^p\to\R^{q}$ by generalising the approach in
 \cite{sawyer:planecurves}.  The plan is that $\psi$ will be a series
 similar to the factorial expansion of $y$, and we hope to make
 $f(y,\psi(y),t)$ close to that value of $f$ where the series for both
 of the first two arguments are truncated---a finite set of values.
 We choose the coefficients in the series to get rid of the main error
 term; it turns out that the coefficients therefore must correspond to
 the values of $\frac{\partial f}{\partial \omega}^{-1}\frac{\partial
 f}{\partial y}$. (The inverse here means a right inverse of the $(n-d)\times q$ matrix $\frac{\partial f}{\partial \omega}$.) So we need to devise a sequence of $q\times p$ real
 matrices that is in some sense `dense' and takes on arbitrarily
 large values, but does not grow too quickly.  This is what we shall
 do now.

 For $k\geq3$ set $$D_k=\left\{(y_2,\dots,y_{k-1})\in
 {}^{(k-2)}\Z^p :0\leq y^i_n\leq n-1\right\},$$ that is, a set of $(k-2)$-tuples of those $p$-dimensional
 vectors that can form the first $k-2$ coefficients in a factorial
 expansion.  Let $\Omega_k$ be the set of all maps
 $D_k\to{}^{q}[-\log\log k,\log\log k]^p$, that is, $q\times
 p$ matrices whose elements are bounded by $\log\log k$.  Next let
 $\{s_j^k\}_{j=1}^{m_k}$ be a finite $1/k$-dense subset of $\Omega_k$,
 meaning that
 $$\forall_{s\in\Omega_k}\exists_j\forall_{(y_2,\dots,y_{k-1})\in
D_k}\left\|s(y_2,\dots,y_{k-1})-s_j^k(y_2,\dots,y_{k-1})\right\|<\frac
1k.$$ At this point it will be helpful to notice that $m_k\sim
(k\log\log k)^{pq(k-1)!^p}$, by taking the number of possible
matrices and raising it to the power of the number of arguments in the
function.

Next we define the sequence of maps to use as coefficients in the
 definition of $\psi$.  Call $r\in\Omega_l$ an {\em extension} of
 $s\in\Omega_k$ if $l\geq k$ and for all $(y_2,\dots,y_{l-1})\in
 D_l$ we have $r(y_2,\dots,y_{l-1})=s(y_2,\dots,y_{k-1})$.
 Set $r_2\equiv 1$ and for each $n\geq3$ choose $r_n\in\Omega_n$ so
 that for all $k\geq3$ and $1\leq j\leq m_k$ there is an $r_n$ that
 is an extension of $s_j^k$.

 Now for $y\in(0,1]^p$ define $\psi$ by
 $$\psi(y)=\sum_{n=2}^\infty r_n(y_2,\dots,y_{n-1})\frac{y_n}{n!}$$
 where  each summand contains a matrix multiplication.  Finally, extend $\psi$  to all of $\R^p$ by periodicity.

 We observe some continuity properties of $\psi$.
\begin{lem}
\label{agree}
 Suppose that $y$ and $\bar{y}$ have the same factorial expansion up
 to the $N$th term (meaning that $y_n=\bar{y}_n$ for
$2\le n\le N$). Then
 $|y-\bar{y}|\leq 1/N!$ and $|\psi(y) -\psi(\bar{y})|\leq
 C\log\log N/N!$.
\end{lem}
\noindent{{\bf Proof: }} $|y-\bar{y}|\leq \sum_{n=N+1}^\infty
\frac{|y_n - \bar{y}_n|}{n!} \leq \sum_{n=N+1}^\infty
\frac{n-1}{n!}=1/N!$. Similarly
\begin{align*}
|\psi(y) -\psi(\bar{y})|&\leq \sum_{n=N+1}^\infty
 \frac{(n-1)\log\log n}{n!}\\ 
&\leq\frac{N\log\log
 (N+1)}{(N+1)!}+\sum_{n=N+2}^\infty \frac{(n-1)\log\log n}{n!}\\ 
&\leq \frac{C\log\log N}{N!}+C\sum_{n=N+2}^\infty
 \frac{n-2}{(n-1)!}\\ 
& = \frac{C\log\log N}{N!}\tag*{\qed}
\end{align*}
 In particular, this shows that $\psi$ is continuous except at points
 where one of the components can also have a terminating factorial
 expansion.  At such points there is left continuity in the ``bad
 components'' and the right limits exist.  Also, $\psi((0,1]^p)$ is a bounded set.
\section{Slices have measure zero}
 We now need to show that for suitable values of $n,d,p,q$ this $\psi$
 has the property claimed, that is, the set
 $$E_f:=\bigcup_y \Gamma \bigl(y,\psi(y)\bigr)$$ has measure zero.  In this
section we show that almost all of the slices through the set at fixed
$t$ have measure zero; since $t$ is fixed we suppress it and just
prove the following:
\begin{lem}
\label{lem:slice}
 Let $f:\R^p\times\R^q\to\R^{n-d}$, with $p\leq n-d\leq q$ and $d<n$. Then if $f$ is $C^1$ and $\frac{\partial f}{\partial \omega}$
 always has highest possible rank (namely $n-d$) and is Lipschitz, then the range of $f(\cdot\,,\psi(\cdot))$ is of
 measure zero.
\end{lem}
 These hypotheses are very natural: $d<n$ is merely to avoid trying to
 pack $n$-dimensional objects in $\R^n$, and the other inequalities
 mean that we should not try to include too large a family of
 surfaces, and we must be free to choose many of the parameters.  The
 condition about the rank simply says that the surface must actually
 depend on the parameters that we are free to vary.
\medskip\par 
\begin{pf}
 By periodicity it is enough to consider only the image of $(0,1]^p$.  For a
 vector $y$ and natural number $k\geq3$ write \begin{align*} y^{(k)}
 &=\sum_{n=2}^{k-1}\frac{y_n}{n!} & \psi^{(k)}(y)&=\sum_{n=2}^{k-1}
 r_n(y_2,\dots,y_{n-1})\frac{y_n}{n!}.  \end{align*} 
Then for natural numbers $k$ and $N$ we have 
\begin{align*}
 f\big(y,&\psi(y)\big)  =
 \left[f\big(y,\psi(y)\big)-f\big(y^{(N)},\psi(y)\big)-\frac{\partial
 f}{\partial
 y}\big(y^{(k)},\psi^{(k)}(y)\big)\big(y-y^{(N)}\big)\right]\\ & +
 \left[f\big(y^{(N)},\psi(y)\big)-f\big(y^{(N)},\psi^{(N)}(y)\big)-\frac{\partial
 f}{\partial
 \omega}\big(y^{(k)},\psi^{(k)}(y)\big)\big(\psi(y)-\psi^{(N)}(y)\big)\right]\\
 & + \left[\frac{\partial f}{\partial
 y}\big(y^{(k)},\psi^{(k)}(y)\big)+\frac{\partial f}{\partial
 \omega}\big(y^{(k)},\psi^{(k)}(y)\big)r_N(y_2,\ldots,y_{N-1})\right]\frac{y_N}{N!}\\
 & + \sum_{n=N+1}^\infty\left[\frac{\partial f}{\partial
 y}\big(y^{(k)},\psi^{(k)}(y)\big)+\frac{\partial f}{\partial
 \omega}\big(y^{(k)},\psi^{(k)}(y)\big)r_n(y_2,\ldots,y_{N-1})\right]\frac{y_n}{n!}\\
 &+ f\big(y^{(N)},\psi^{(N)}(y)\big)\\
 &=:I(y)+II(y)+III(y)+IV(y)+V(y).  \end{align*} The final term
 takes a very large, but finite, number of values, so our task is to
 show that the other terms are correspondingly extremely small.

Let $\eps>0$ be given.  Using the hypothesis that $f$ is $C^1$ together with the fact that $y$ and $\omega$ lie in the bounded sets $(0,1]^p$ and $\psi((0,1]^p)$ respectively, choose
$k$ so large that the following hold:
\begin{enumerate}\label{k-enum}
 \item \label{iy}If both $|y-\bar{y}|<\frac{1}{(k-1)!}$ and
 $|\omega-\bar{\omega}|<\frac{\log\log k}{(k-1)!}$, then 
 $\left\|\frac{\partial f}{\partial y}(y,\omega)-\frac{\partial
 f}{\partial y}(\bar{y},\bar{\omega})\right\| < \eps$.  This is possible since the Jacobian is continuous.
 \item \label{iw}If both $|y-\bar{y}|<\frac{1}{(k-1)!}$ and
 $|\omega-\bar{\omega}|<\frac{\log\log k}{(k-1)!}$, then 
 $\Bigl\|\frac{\partial
 f}{\partial \omega}(y,\omega)-\frac{\partial f}{\partial
 \omega}(\bar{y},\bar{\omega})\Bigr\|<\eps/(k\log k)$. This is possible since the Jacobian is Lipschitz.  
\item \label{ii}$\left\|\frac{\partial f}{\partial
y}(y,\omega)\right\|<\log\log k$ and $\Bigl\|\frac{\partial
f}{\partial \omega}(y,\omega)\Bigr\|<\log\log k$.
\item $\left\|\frac{\partial f}{\partial
\omega}^{-1}\frac{\partial f}{\partial
y}(y,\omega)\right\|<\log\log k$ where $\frac{\partial f}{\partial
\omega}^{-1}$ is a right inverse of the $(n-d)\times q$ matrix
$\frac{\partial f}{\partial \omega}$.  Here we are using the
assumptions that $q\geq n-d$ and that the matrix has full rank.
\item \label{iv}$\frac{(\log\log k)^2}{k}<\eps$ \end{enumerate}
Next, find an $s_j^k$ within $1/k$ of the matrix $-\frac{\partial
f}{\partial \omega}^{-1}\frac{\partial f}{\partial
y}(y^{(k)},\psi^{(k)}(y))$.  Then find $N$ such that $r_N$ is an
extension of $s_j^k$.  We show that parts $I$--$IV$ above are smaller
than $\frac{C\eps}{(N-1)!}$.

 Part $I$ is handled using the mean value theorem.
 The $i$th component of $I(y)$ is 
 $$f^i\bigl(y,\psi(y)\bigr)-f^i\bigl(y^{(N)},\psi(y)\bigr)-\grad_y f^i\bigl(y^{(k)},\psi^{(k)}(y)\bigr)\cdot(y-y^{(N)})$$
  which, by the one-dimensional mean value theorem in the direction $y-y^{(N)}$, equals 
$$\left(\grad f^i\bigl(\xi,\psi(y)\bigr)-\grad f^i\bigl(y^{(k)},\psi^{(k)}(y)\bigr)\right)\cdot(y-y^{(N)})$$ 
for some $\xi\in[y^{(N)},y]$.  But then $|\xi-y^{(k)}|<\frac{1}{(k-1)!}$, and $|\psi(y)-\psi^{(k)}(y)|<\frac{\log\log k}{(k-1)!}$ so that by applying \ref{iy} to this and all the other components we eventually get
 \begin{align*}
 |I(y)| & \leq \eps |y-y^{(N)}|\\
 & \leq \eps \sum_N^\infty\frac{1}{n!}|y_n|\\
 & \leq \eps \frac{1}{(N-1)!}.
 \end{align*}
$II$ works similarly, except that we end up with 
 $$|II(y)|\leq\frac\eps{k\log k}|\psi(y)-\psi^{(N)}(y)|\leq C\frac{\eps}{k\log k}\frac{\log\log
N}{(N-1)!}.$$ 
But note that $N$ was chosen to make $r_N$ an extension of $s_j^k$, so
that provided we ordered the sequence $(r_n)$ sensibly, we have
\begin{align*} N & \leq \sum_{l<k}m_l+j\\ & \leq Ck(k\log\log
k)^{pq(k-1)!^p} \end{align*} 
and hence
$\log\log N\lesssim k\log k$.  So the estimate of
$C\frac{\eps}{(N-1)!}$ for $II$ follows.  (This is the step for which we need the rather unlikely-looking double log---in Sawyer's proof this
issue does not arise, because there $\frac{\partial f}{\partial
\omega}$ is minus the identity matrix and so this whole term is zero.)

 For $III$, our choice of $N$ gives us cancellation.
 \begin{align*}
 |III(y)| & \leq 0+\frac{1}{k}\left\|\frac{\partial f}{\partial
 \omega}\bigl(y^{(k)},\psi^{(k)}(y)\bigr)\right\|\frac{|y_N|}{N!}\\ & \leq
 \frac{\log\log k}{k}\frac{N-1}{N!}\\ & < \frac{\eps}{(N-1)!}.
 \end{align*}
 Finally,
 \begin{align*}
 IV(y)& \leq\sum_{n=N+1}^\infty \frac{\log\log k (1+\log\log n)(n-1)}{n!}\\
 &\leq \log\log k\left[\frac{CN\log\log N}{(N+1)!}+\sum_{N+2}\frac{(n-1)\log\log n}{n!}\right]\\
 &\leq \frac{C\log\log k \log\log N}{N!}\\
 &\leq \frac{C(\log\log k)^2}{k(N-1)!}\qquad\mbox{since $N>k$}\\
 &< \frac{C\eps}{(N-1)!}.
 \end{align*}
 Combining these estimates we see that
 $$\operatorname{range}\big(f(\cdot\,,\psi(\cdot))\big)\subseteq\bigcup_{z\in\operatorname{range}(V)}B\left(z,{\tfrac{C\eps}{(N-1)!}}\right).$$ But $V(y)$ depends only on
$y_2,\dots,y_{N-1}$, so $\operatorname{range}(V)$ has at most
$(N-1)!^p$ elements.  Hence 
\begin{align*}
\left|\operatorname{range}\big(f(\cdot\,,\psi(\cdot))\big)\right| & \leq
(N-1)!^p\left(\frac{C\eps}{(N-1)!}\right)^{n-d}\\ & =
C\frac{\eps^{n-d}}{(N-1)!^{n-d-p}} \end{align*} 
which,
since $\eps$ is arbitrary, proves the result since $p\leq n-d$ and
$d<n$.  \end{pf}
{\noindent{\bf Proof of Theorem~\ref{thm:null}: }}
 To conclude the proof of the theorem we must show that the entire set
 $E_f$ is measurable.  
Consider the set 
$$E:=\bigcup_{(y_2,y_3,\dots)}\bigcap_{k=3}^{\infty}\left\{\begin{pmatrix}x\\t\end{pmatrix}:\left|x-f\bigl(y^{(k)},\psi^{(k)}(y),t\bigr)\right|\leq\frac1k\right\}$$
where the union is taken over all infinite sequences $(y_2,y_3,y_4,\dots )$ with each vector $y_m$ belonging to $\{0,1,\dots,m-1\}^p$. The sets intersected are measurable sets depending only on the first $k-2$ terms of the sequence; therefore $E$ is the result of applying the Souslin operation to a class of measurable sets and hence (see for example \cite[page~45]{rogers:hausdorff}) is measurable.
Since $(y,\omega)\mapsto f(y,\omega,t)$ is $C^1$ for a.e.~$t$, the set $E$ is just the union of the surfaces $\Gamma\bigl(y,\psi(y)\bigr)$ except at those $t$ for which $f$ is not $C^1$.  That is, $E$ differs from $E_f$ only on a set of measure zero.  Therefore $E_f$ is measurable.{\qed\medskip\par}
 \section{Discussion}
We remark that our hypotheses are stronger
 than needed: The Lipschitz condition on $\frac{\partial f}{\partial\omega}$ was only used to show that given $\eps>0$ we can find $k$ such that \ref{iw} is true:
 this would still hold with a weaker condition on the modulus of
 continuity of the Jacobian. Moreover, by replacing $\log\log$ throughout the proof by three or more logs, we could weaken the condition further. 
In fact, we could do without any such
 condition if we sacrificed the universality of $\psi$ and allowed it to
 depend on the rate of growth of the derivatives of $f$.  It may also be possible to relax the $C^1$ hypothesis slightly, although Sawyer shows that it cannot be replaced by a Lipschitz condition of any order less than $1$.

 Our theorem sheds some light on
 other known results on curve-packing.  For example, a null set in the
 plane can be constructed so as to include a circle of every radius
 (Besicovitch and Rado, Kinney 1968), but if a set has a circle
 centred at every point in the plane then it must have positive
 measure (Bourgain 1986, Marstrand 1987).  However, with circles
 centred at all points on a curve the set can still be null (Talagrand
 1980).  These examples illustrate the numerology of the theorem and
 suggest that the conditions on the parameters might in fact be
 necessary as well as sufficient.  Higher dimensional examples include
 the $k$-plane problem: A set in $\R^3$ that includes a plane in
 every direction must have positive measure (Marstrand 1979, Falconer
 1980)---what can be said about packing $k$-planes in $\R^n$?  This
 problem has been studied by Falconer, Bourgain and others but remains
 unsolved.  In this case we would have $d=k,\ p=k(n-k)$ and $q=n-k$, so that if the numerology of Theorem~\ref{thm:null} was found to
 be sharp, then $k$-planes could be packed into a null set only when
 $k=1$. (Towards this, Mitsis \cite{mitsis:2plane} has recently shown that a set in $\R^n$ containing a translate of every $2$-plane must have full dimension.) References for these and similar results can be found in
 \cite{falconer:frac}, \cite{mattila:geomsets} and
 \cite{wolff:recentkakeya}.
\paragraph{Acknowledgements} This work formed part of my PhD thesis \cite{wisewell:thesis}, which also contains results on the dimension of curved Kakeya and Nikodym sets as well as relating these to open problems on oscillatory integrals.  The encouragement and help of my supervisor Professor A.~Carbery, and the financial support of the EPSRC and the Seggie-Brown Trust are gratefully acknowledged.
%\bibliography{/home/laura/latex/all}

\begin{thebibliography}{1}

\bibitem{falconer:frac}
K.~J. Falconer.
\newblock {\em The geometry of fractal sets\/}.
\newblock No.~85 in Cambridge Tracts in Mathematics. Cambridge University
  Press, 1985.

\bibitem{mattila:geomsets}
P.~Mattila.
\newblock {\em Geometry of Sets and Measures in {E}uclidean Spaces\/}.
\newblock No.~44 in Studies in Advanced Mathematics. Cambridge, 1995.

\bibitem{mitsis:2plane}
T.~Mitsis.
\newblock $(n,2)$-sets have full {H}ausdorff dimension.
\newblock {\em To appear in Revista Matem{\'a}tica Iberoamericana\/}, 2003.

\bibitem{rogers:hausdorff}
C.~A. Rogers.
\newblock {\em Hausdorff measures\/}.
\newblock Cambridge University Press, Cambridge, 1998.
\newblock Reprint of the 1970 original, with a foreword by K. J. Falconer.

\bibitem{sawyer:planecurves}
E.~Sawyer.
\newblock Families of plane curves having translates in a set of measure zero.
\newblock {\em Mathematika\/}, {\bf 34}: 69--76, 1987.

\bibitem{wisewell:thesis}
L.~Wisewell.
\newblock {\em Oscillatory Integrals and Curved {K}akeya Sets\/}.
\newblock Ph.D. thesis, University of Edinburgh, 2003.

\bibitem{wolff:recentkakeya}
T.~Wolff.
\newblock Recent work connected with the {K}akeya problem.
\newblock In H.~Rossi, ed., {\em Prospects in Mathematics (Princeton, New
  Jersey, 1996)\/}, pp. 129--162. American Mathematical Society, Providence,
  RI, 1999.

\end{thebibliography}
%\bibliographystyle{lauraabbrv2}

\bigskip
\noindent Laura Wisewell\\
School of Mathematics\\
University of Edinburgh\\
King's Buildings\\
Edinburgh\\
EH9 3JZ\\
Scotland\\
United Kingdom\smallskip\par
\noindent {\tt l.wisewell@ed.ac.uk}

\end{document}